\documentclass{article}
\usepackage{graphicx} 
\usepackage{tikz}
\usepackage{amssymb}
\usepackage{amsmath}
\usepackage{enumitem}
\usepackage{framed}
\usepackage{xcolor}
\usepackage{listings}
\usepackage{url}
\usepackage{authblk}

\definecolor{keywordcolor}{rgb}{0,0,0}
\definecolor{tacticcolor}{rgb}{0,0,0}
\definecolor{symbolcolor}{rgb}{0,0,0}
\definecolor{sortcolor}{rgb}{0,0,0}
\definecolor{attributecolor}{rgb}{0,0,0}

\title{Capturing properties of planar diagrams in Lean proof assistant software}
\author[1]{Alastair Litterick} 
\author[2]{Alexei Vernitski} 
\author[3]{Billy Woods}
\affil[1]{Mathematics, Statistics and Actuarial Sciences, University of Essex}
\affil[2]{Mathematics, Statistics and Actuarial Sciences, University of Essex. Corresponding author asvern@essex.ac.uk}
\affil[3]{Essex Pathways, University of Essex}
\begin{document}
\maketitle
\begin{abstract}
Automated proof assistants are a technology pre-empting mistakes in mathematics. In our practice we have seen that reasoning about planar diagrams is difficult to both humans and computers. One example that has led to wrong statements in publications is that an orientation-preserving mapping is not always defined by how it acts on triples of elements. In this paper we formalise orientation-preserving mappings in proof assistant software Lean and report on our take-aways.
\end{abstract}


\section{Introduction}

In our practice we have observed that making conclusions concerning planar diagrams is tricky and frequently leads to mistakes. Famous examples include the controversy regarding the first proof of the Jordan curve theorem and Kempe's incorrect proof of the four-colour theorem. Examples known to us from our own recent practice include \cite{lisitsa2023describing} and \cite{howie2020planar} correcting errors in earlier papers of other authors. 

The example considered in this article is also taken from our recent practice. Order-preserving mappings (which is how monotonic functions are called in semigroup theory) are defined by how they act on pairs of elements of the chain, that is, a mapping $f$ is called order-preserving if from $a \le b$ it follows that $f(a) \le f(b)$. Orientation-preserving mappings are a generalization of order-preserving mappings. As we will see in  the following section, orientation-preserving mappings are `almost' defined by how they act on triples of elements of the chain, but not quite; there are exceptions. Relevant publications where this distinction is studied are \cite{levimitchell2006,higgins2022orientation,fernandes2025groups}. Although it is well known that orientation-preserving mappings are `not quite' defined by how they act on triples of elements of the chain, this fact frequently gets overlooked and has led to small mistakes in some published results; indeed, paper \cite{higgins2022orientation} required a corrigendum \cite{higgins2024correction}, and paper \cite{fernandes2023CKMS} required a corrigendum \cite{fernandes2024CKMS}.

Automated proof assistants are a kind of technology that ensures that proofs exactly match the intended results, thus pre-empting mistakes. At the same time, it is difficult to write code for an automated proof assistant that can be used to reason about planar diagrams. Are automated proof assistants sufficiently developed to verify facts like the ones described above? 
In this paper we formalise orientation-preserving mappings in proof assistant software Lean and report on our take-aways.

Note that within this article, it is not our purpose to define planarity as such. For our purposes, planarity is a sweet spot where both humans and computers, for different reasons, find it difficult to make correct inferences, and therefore, this is where we can draw examples showing how human reasoning and automated proof assistants can complement one another.

\section{Orientation-preserving mappings}\label{sec:mappings}

Fix a positive integer $n$. The set $\{0,1,\dots,n-1\}$ is denoted by $[n]$. By a mapping we mean a transformation on $[n]$. For simpler notation, a mapping $0\mapsto a_0,1\mapsto a_1,\dots,n-1\mapsto a_{n-1}$ will be denoted by $(a_0,a_1,\dots,a_{n-1})$.

The standard definition of an orientation-preserving mapping is as follows, see, for example, \cite{higgins2022orientation, fernandes2023CKMS,fernandes2025groups}. Let $s=(a_1,a_2,\dots,a_t)$ be a sequence of $t$ elements of $[n]$. We say that $s$ is cyclic [respectively, anti-cyclic] if there is no more than one index $i\in\{1,\dots,t\}$ such that $a_i>a_{i+1}$ [respectively, $a_i<a_{i+1}$], where $a_{t+1}$ denotes $a_1$.  

Given a mapping $f$ whose domain is $\{a_1<\dots<a_t \}$, we say that $f$ is orientation-preserving (respectively, orientation-reversing) if the sequence of its images $(f(a_1), f(a_2),\dots, f(a_t))$ is cyclic (respectively, anti-cyclic).

Now we will comment on the definition of an orientation-preserving mappings, concentrating on intuitive understanding. In the following sections, when we formalise this concept for the automated proof assistant Lean, we lean in the opposite direction, making the concept as formal as it can be.

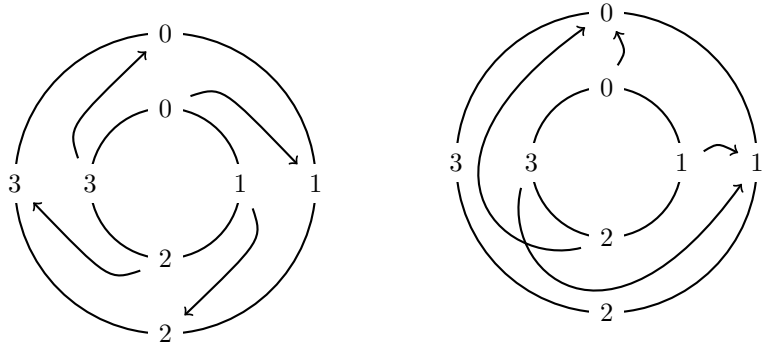
\begin{figure}
\centering 
\begin{minipage}{.45\textwidth}
\centering 

\begin{tikzpicture}[scale=1]

\def\R{2}   
\def\r{1}   

\draw[thick] (0,0) circle (\R);
\draw[thick] (0,0) circle (\r);

\foreach \i/\a in {0/90, 1/0, 2/270, 3/180} {
  \draw[thick,->,shorten >=10pt,shorten <=10pt]
    (\a:\r+0.05) .. controls (\a-30:{(\R+\r)/2}) .. (\a-90:\R);
    
  \node[fill=white, inner sep=4pt] at (\a:\r) {\i};
  \node[fill=white, inner sep=4pt] at (\a:\R) {\i};
}

\end{tikzpicture}
\end{minipage}
\begin{minipage}{.45\textwidth}
\centering 
\begin{tikzpicture}[scale=1]

\def\R{2}   
\def\r{1}   

\draw[thick] (0,0) circle (\R);
\draw[thick] (0,0) circle (\r);

\foreach \a/\b/\c/\d in {90/80/80/90, 0/10/10/0} {
  \draw[thick,->,shorten >=8pt,shorten <=8pt]
    (\a:\r+0.05) .. controls (\b:{(\R+\r)/2}) and (\c:{(\R+\r)/2}) .. (\d:\R);
}    
\foreach \a/\b/\c/\d in {270/225/180/90, 180/225/270/0} {
  \draw[thick,->,shorten >=10pt,shorten <=10pt]
    (\a:\r+0.05) .. controls (\b:{(\R+\r)/2+0.5}) and (\c:{(\R+\r)/2+1.3}) .. (\d:\R);
}    
\foreach \i/\a in {0/90, 1/0, 2/270, 3/180} {
  \node[fill=white, inner sep=4pt] at (\a:\r) {\i};
  \node[fill=white, inner sep=4pt] at (\a:\R) {\i};
}
\end{tikzpicture}
\end{minipage}
\caption{The mapping on the left is orientation-preserving. The mapping on the right is not orientation-preserving.}
\label{fig:examples}
\end{figure}

Draw $[n]$ around a circle clockwise. Around that circle, draw another copy of $[n]$ around a circle, also clockwise. Then draw connectors induced by the mapping. If the connectors can be drawn in such a way that they do not intersect, the mapping is orientation-preserving; for example, mapping $(1,2,3,0)$ shown in Figure \ref{fig:examples} on the left, is orientation-preserving. If the connectors intersect in whichever way you lay them out, the mapping is not orientation-preserving; for example, mapping $(0,1,0,1)$ shown in Figure \ref{fig:examples} on the right, is not orientation-preserving. Mapping $(0,1,0,1)$ is an example demonstrating that an orientation-preserving mapping is not always defined by how it acts on triples of elements; indeed, every restriction of mapping $(0,1,0,1)$ to $3$ elements is orientation-preserving. As we discussed in the previous section, this counterexample, although well-known, was overlooked in the statements of some results in \cite{higgins2022orientation, fernandes2023CKMS}.

To describe the concept of an orientation-reversing mapping, draw $[n]$ around a circle clockwise, and around that circle, draw another copy of $[n]$ around a circle, but anti-clockwise. Then draw connectors induced by the mapping. If the connectors can be drawn in such a way that they do not intersect, the mapping is orientation-reversing. If the connectors intersect in whichever way you lay them out, the mapping is not orientation-reversing.

In the following sections we use proof assistant software Lean to demonstrate that sequence $(0,1,0,1)$ is neither cyclic nor anti-cyclic (we do not expand this result to proving that mapping $(0,1,0,1)$ is neither orientation-preserving nor orientation-reversing because this would distract us from our main topic of making inferences about planar diagrams). We then proceed to considering a theorem \cite[Theorem 3]{higgins2022orientation} stating that if a sequence does not include exactly $2$ distinct elements then its being cyclic (or anti-cyclic) is determined by the orientations of its subsequences of length $3$.


\section{Lean}

\subsection{What are Lean and mathlib?}

Lean 4 \cite{lean4} is the current major version of the Lean programming language. In the context of mathematical research, Lean is often referred to as a \emph{proof assistant} or an \emph{interactive theorem prover} (though these terms may give an inflated impression of its current capabilities in the \emph{production} of proofs). The primary current use of Lean in this context is in proof \emph{verification}: the researcher will attempt to formalise a theorem in Lean, i.e.\ provide a fully elaborated proof with completely explicit, checkable steps which can be formally verified by Lean's typechecking logic. If the verification fails, Lean will attempt to provide useful error messages for the researcher.

mathlib4 is the current version of mathlib, a project maintained by the Lean community \cite{mathlib}, and is a library of formalised mathematics built in Lean. It contains a very large base of definitions and theorems, across most areas of undergraduate mathematics, which the researcher can build on and adapt. As a well-written repository, it is also an excellent (if somewhat advanced) resource for learning how to formalise theorems in Lean; the more beginner-friendly interactive online tutorial book Mathematics in Lean \cite{mil} is built on top of mathlib.

The distinction is (roughly) that Lean's core and standard library contain what is relevant to programming, and mathlib contains what is relevant to mathematics. Lean provides the basic logical language (logical operators and quantifiers), elementary data structures (natural numbers, lists, booleans), some basic operations on these objects (arithmetic on natural numbers, appending or rotating lists), and some basic ``proof tactics" (simplification using basic rules, rewriting an expression using a known equality, induction). In contrast, mathlib contains notions such as groups and topological spaces, as well as mathematical results.

The distinction can be seen clearly in some of mathlib's basic theorems. As an example: natural numbers (\texttt{Nat}) are defined within Lean, but the notion of a natural number being prime (\texttt{Nat.Prime}) is only defined within mathlib, and so the theorem that 2 is a prime number (\texttt{Nat.prime\_two}) is contained within mathlib. As another example: integers (\texttt{Int}) and their basic arithmetic (e.g.\ multiplication) are defined within Lean, but the theorem that multiplication is commutative (\texttt{Int.mul\_comm}), or that multiplication by $-1$ is the same as negation (\texttt{Int.mul\_neg\_one}), are contained within mathlib.

\subsection{Writing a proof in Lean}

As the above examples suggest, writing a proof in Lean can involve significant overhead, as even the most ``obvious" mathematical steps must be explicitly encoded. In particular, arguments which rely on intuition (such as many geometric arguments) cannot be encoded, and many simple proofs (which rely on the reader to fill in the gaps) become long and tedious to formalise.

This is both a weakness and a strength of formalised mathematics. It can be very time-consuming to formalise a result in Lean, and many researchers may decide that the time required to successfully formalise a result exceeds the benefit. On the other hand, the omission of these ``obvious" steps in research papers can lead to very subtle mistakes scattered throughout the literature, which may themselves take years or decades for the mathematical community to notice and correct. The cost-benefit tradeoff is currently unclear, but there are clear arguments in favour of Lean's very high standard for guaranteed correctness.

We mention a few points of friction surrounding automation. We have found Lean poorly suited to ``exploratory" or ``interactive" reasoning: arguments must first be spelt out in detail on paper, then translated into Lean syntax. The automation available consists of a large collection of separate tools, of which we mention a few below:

\begin{itemize}
\item The built-in simplifier \lstinline{simp} recursively rewrites expressions based on a series of \emph{rewrite rules}, aimed at putting everything into a syntactic normal form. The rewrite rules include a carefully curated set of lemmas, which aim to rewrite ``more complex" expressions as ``simpler" expressions. The researcher can extend \lstinline{simp} by adding extra lemmas, but this is delicate: lemmas must be chosen and written carefully, so as to avoid unwanted results such as infinite loops.

\lstinline{simp} can be rather a blunt instrument, and the normal form it reaches may not always be what is desired: in practice, we have found this rather unpredictable. \lstinline{simp} functions as a black box, and does not report exactly \emph{how} it simplified an expression, only giving the final simplified form. There are many variants, such as \lstinline{simp only} (which allows the user to restrict to certain simplification rules only) and \lstinline{simp?} (which locates potential rewrites and then suggests them to the user), which overcome this issue at the expense of automation.
\item The pattern-matching search \lstinline{apply?} will attempt to provide suggestions for theorems that can be applied. This appears to be most useful when the user has forgotten the name of a lemma.
\item There are a number of domain-specific solvers, such as \lstinline{linarith}, a complete decision procedure which can solve many goals in linear arithmetic using algorithms such as Fourier-Motzkin elimination. This is very powerful, but only within a very narrow scope.
\item The more general tool \lstinline{aesop} aims to prove any goal by searching through the space of possible proofs using complex search heuristics. It is particularly good at structural manipulation, such as propositional logic, proofs with multiple cases, and so on.
\end{itemize}

Unfortunately, the overall result is not currently powerful enough to perform sophisticated reasoning automatically: in our experience, the primary use of these tools has been in overcoming barriers inherent to Lean, not in overcoming mathematical barriers. Moreover, while these tactics are powerful in their own domains, there does not seem to be any useful unification of them so far. It is unclear how (if at all) this will change in the future: the recent \lstinline{grind} tactic, which integrates techniques motivated by modern SMT solvers, shows promise, but is still in the early stages of development.

\subsection{Inductive types}

Like many other proof assistants, Lean is built around the \emph{calculus of constructions}, a type theory. A full discussion of type theory is beyond the scope of this note, but we mention a few aspects.

Every expression in Lean has an associated \emph{type}. Many of these will be familiar:
\begin{itemize}[noitemsep]
\item \lstinline{3} has type \lstinline{Nat} (natural number),
\item \lstinline{false} has type \lstinline{Bool} (boolean),
\item \lstinline{3 = 4} has type \lstinline{Prop} (proposition).
\end{itemize}
There are many ways of defining new types from old: for instance, using the Lean operator \lstinline{→}, we may encode functions from $\mathbb{N}$ to $\mathbb{N}$ as expressions of type \lstinline{Nat → Nat}.

\emph{Inductive} types are built from a list of constructors. For instance, \lstinline{Nat} is an inductive type, and terms of type \lstinline{Nat} may be built using one of two constructors:
\begin{itemize}[noitemsep]
\item \lstinline{zero} (representing 0),
\item \lstinline{succ}, which takes an additional parameter \lstinline{n} of type \lstinline{Nat} (representing the successor of \lstinline{n}).
\end{itemize}

Theorems about a natural number \lstinline{n} are then proven inductively by splitting them into cases: the case where \lstinline{n} is \lstinline{zero}, and the case where \lstinline{n} is \lstinline{succ m}.

\subsection{The code: $(0,1,0,1)$ is neither cyclic nor anti-cyclic}

The following represents our attempt to formalise this small result in Lean 4: we do not claim that this is the most efficient or idiomatic way to do so. (In fact, we will see in \S \ref{subsec:nextsteps} below that it is not well-suited to all purposes.) All code below was developed and tested in Lean 4.27.0 with mathlib 4.27.0.

We must first teach Lean to reason about sequences in the sense of \S \ref{sec:mappings}, which we encode as terms of type \lstinline{List Nat} (lists of natural numbers). Like the type \lstinline{Nat}, the type \lstinline{List Nat} is defined inductively: its terms may be either
\begin{itemize}[noitemsep]
\item \lstinline{nil}, representing an empty list of natural numbers, or
\item \lstinline{cons a as} (where \lstinline{a} has type \lstinline{Nat} and \lstinline{as} has type \lstinline{List Nat}), representing a concatenation: we may think of this as
$(\text{\lstinline{a}}, \underbrace{x_2, \dots, x_n}_{\text{\lstinline{as}}}).$
\end{itemize}
Hence the natural way to work with this type is to ``peel off" one element at a time from the left. This makes it easy to define what it means for a sequence to be \emph{increasing}:

\begin{framed}
\begin{lstlisting}[language=lean]
def List.isIncreasing : List Nat → Bool
  | .nil                 => true
  | .cons _ .nil         => true
  | .cons a (.cons b bs) => a ≤ b && (cons b bs).isIncreasing
\end{lstlisting}
\end{framed}

In other words: \lstinline{.isIncreasing} is a function that takes in a sequence and returns a boolean value, which is determined inductively by the following:
\begin{itemize}[noitemsep]
\item the empty list $()$ returns \lstinline{true} (i.e.\ the empty sequence is increasing),
\item the concatenation of any natural number and the empty list returns \lstinline{true} (i.e.\ a sequence containing only one natural number is increasing),
\item $(\text{\lstinline{a}}, \text{\lstinline{b}}, x_3, \dots, x_n)$ returns \lstinline{true} if and only if both \lstinline{a ≤ b} and $(\text{\lstinline{b}}, x_3, \dots, x_n)$ is increasing.
\end{itemize}

To define the term \emph{cyclic}, we must have access to both the first and the last elements of the sequence at the same time. For instance, we may \emph{augment} the sequence $(x_1, x_2, \dots, x_n)$ to the sequence $(x_n, x_1, x_2, \dots, x_n)$ as follows.

\begin{framed}
\begin{lstlisting}[language=lean]
def augmentList : List Nat → List Nat
  | .nil         => .nil
  | .cons a as   =>
    if h : as.reverse.length > 0 then
      .cons (as.reverse.get ⟨0, h⟩) (.cons a as) else
      .cons a .nil
\end{lstlisting}
\end{framed}

In other words:
\begin{itemize}[noitemsep]
\item augmenting the empty sequence returns the empty sequence,
\item otherwise, to augment $(\text{\lstinline{a}}, \underbrace{\dots}_{\text{\lstinline{as}}})$, first check whether \lstinline{as} is empty.
\begin{itemize}
\item If not, adjoin the final element of \lstinline{as} to $(\text{\lstinline{a}}, \underbrace{\dots}_{\text{\lstinline{as}}})$ and return the result.
\item If so, just return $(\text{\lstinline{a}})$.
\end{itemize}
\end{itemize}

Here we see an example of the precision required by Lean. On line 5 of this snippet, to obtain the \emph{final} element of the sequence \lstinline{as}, we first reverse it to obtain the new sequence \lstinline{as.reverse}, then attempt to \lstinline{.get} its \emph{first} element. (The indexing begins from 0.) However, we cannot write \lstinline{as.reverse.get 0}: since it is not possible to access the first element of a list of length 0, Lean also requires an assurance that \lstinline{as.reverse} has length greater than 0. If this hypothesis is true, we call it \lstinline{h} (on line 4) and pass it to \lstinline{.get} (on line 5).

Now we can define what it means to be cyclic. We do this via an auxiliary definition:

\begin{framed}
\begin{lstlisting}[language=lean]
def List.isCyclicAux : List Nat → Bool
  | .nil                 => true
  | .cons _ .nil         => true
  | .cons a (.cons b bs) =>
    if a ≤ b then
      (cons b bs).isCyclicAux else
      (cons b bs).isIncreasing
      
def List.isCyclic (S : List Nat) : Bool :=
  (augmentList S).isCyclicAux
\end{lstlisting}
\end{framed}

In other words, to check whether a sequence is \emph{cyclic-aux}, we proceed inductively:
\begin{itemize}[noitemsep]
\item If the sequence is empty, it is cyclic-aux.
\item Otherwise, if the sequence is a singleton, it is cyclic-aux.
\item Otherwise, write the sequence as $(\text{\lstinline{a}}, \text{\lstinline{b}}, x_3, \dots, x_n)$.
\begin{itemize}
\item \raggedright If \lstinline{a ≤ b}, then this sequence is cyclic-aux if and only if the subsequence $(\text{\lstinline{b}}, x_3, \dots, x_n)$ is cyclic-aux.
\item \raggedright If \lstinline{a > b}, then this sequence is cyclic-aux if and only if the subsequence $(\text{\lstinline{b}}, x_3, \dots, x_n)$ is increasing.
\end{itemize}
\end{itemize}
(This encodes the same idea as the original definition: there can only be at most \emph{one} decreasing adjacent pair. Of course, many mathematically equivalent formulations of this definition are possible, and much of the difficulty lies in choosing the definition that is easiest to work with in Lean.)

Finally, a sequence is cyclic if its augmentation is cyclic-aux.

We may similarly define:

\begin{framed}
\begin{lstlisting}[language=lean]
def List.isAnticyclic (S : List Nat) : Bool :=
  (S.reverse).isCyclic
\end{lstlisting}
\end{framed}

Next, we define four possibilities for the \emph{orientation sort} of a sequence:

\begin{framed}
\begin{lstlisting}[language=lean]
inductive orientationSort where
  | none       : orientationSort
  | cyclic     : orientationSort
  | anticyclic : orientationSort
  | both       : orientationSort
\end{lstlisting}
\end{framed}

The following definition allows Lean to parse an explicit sequence and determine which orientation sort it has:

\begin{framed}
\begin{lstlisting}[language=lean]
def List.orientation (S : List Nat) : orientationSort :=
  if S.isCyclic then
    if S.isAnticyclic then
      orientationSort.both else
      orientationSort.cyclic
  else
    if S.isAnticyclic then
      orientationSort.anticyclic else
      orientationSort.none
\end{lstlisting}
\end{framed}

Now we may finally ask Lean to evaluate the orientation of the sequence $(0, 1, 0, 1)$:

\begin{framed}
\begin{lstlisting}[language=lean]
#eval [0, 1, 0, 1].orientation
\end{lstlisting}
\end{framed}

The Lean parser responds with \lstinline{orientationSort.none}, confirming that this sequence is neither cyclic nor anti-cyclic.

\subsection{Next steps: proving a theorem}\label{subsec:nextsteps}

So far, most of our work has been definitions; a sensible next step from here would be to formalise a small theorem, such as \cite[Theorem 3]{higgins2022orientation} (after the corrigendum of \cite{higgins2024correction}), which says that the orientation of a sequence of rank $\neq 2$ is determined by the orientations of its subsequences of length 3. It is usually a good idea to break a theorem down into the smallest possible lemmas: for instance, the user may begin by attempting to prove that a subsequence of a cyclic sequence is cyclic, or even that a subsequence of an increasing sequence is increasing. Here, we use the built-in definition \lstinline{Sublist}.

\begin{framed}
\begin{lstlisting}[language=lean]
theorem sublist_of_increasing_is_increasing
  (S T : List Nat)  (hS : S.isIncreasing) (hTS : T.Sublist S) :
  T.isIncreasing := ...
\end{lstlisting}
\end{framed}

(The terms on line 2 are declarations and hypotheses: ``\emph{let} \lstinline{S} and \lstinline{T} be sequences, and \emph{assume} that \lstinline{S} is increasing, and \emph{assume} that \lstinline{T} is a subsequence of \lstinline{S}". Here, \lstinline{hS} and \lstinline{hTS} are names for the respective hypotheses by which we may refer to them later. The term on line 3 is the goal we wish to prove. The proof will follow \lstinline{:=}.)

The proof is case-checking: the cases where either \lstinline{S} or \lstinline{T} is \lstinline{nil} are easy, and then we must deal with the cases where \lstinline{S = cons x S'} and \lstinline{T = cons y T'}. This then naturally splits further into the two cases \lstinline{x = y} and \lstinline{x ≠ y}.

Unfortunately, any attempt to fill in the details of this proof will take dozens of lines of code, for two reasons. Firstly, our definition of \lstinline{isIncreasing} above outputs a \lstinline{Bool}, which is well-suited to \emph{calculation}, but not to inductive \emph{proofs}. Secondly, as our definition is entirely custom, there are no ready-proved lemmas to draw on.

A much shorter and more idiomatic Lean \emph{proof} might go as follows:

\begin{framed}
\begin{lstlisting}
import Batteries.Data.List.Basic
import Mathlib.Data.List.Chain

def List.isIncreasing (S : List Nat) : Prop :=
  IsChain (· ≤ ·) S

theorem sublist_of_increasing_is_increasing
  (S T : List Nat)  (hS : S.isIncreasing) (hTS : T.Sublist S) :
  T.isIncreasing := by
    simp [List.isIncreasing]
    exact hS.sublist hTS
\end{lstlisting}
\end{framed}

We comment briefly on the code.
\begin{itemize}[noitemsep]
\item We begin by importing two modules, \lstinline{Batteries.Data.List.Basic} and \\ \lstinline{Mathlib.Data.List.Chain}, which allow us to access the pre-defined proposition \lstinline{IsChain} and helper lemma \lstinline{IsChain.sublist} respectively.
\item Our definition now says roughly that a sequence \lstinline{S} is increasing if the relation \lstinline{(· ≤ ·)} holds between adjacent elements of \lstinline{S}.
\item In the theorem, to prove that \lstinline{T} is increasing, we first simplify using the definition of \lstinline{List.isIncreasing}. This expands the definition, and replaces it with the simpler \lstinline{IsChain} definition: our goal then becomes to prove that \lstinline{T} is a chain of elements satisfying \lstinline{(· ≤ ·)}. But this is proved exactly by the general helper lemma \lstinline{IsChain.sublist} applied to our hypotheses.
\end{itemize}

To summarise, on the one hand, in this section we gave an overview of a selection of Lean tools that can be successfully used to make progress towards proving \cite[Theorem 3]{higgins2022orientation}, and on the other hand, we saw that actually proving \cite[Theorem 3]{higgins2022orientation} in Lean would require inordinate amount of work, so we did not proceed to writing complete code for this proof.

\subsection{Adopting Lean in research and teaching}

Lean comes with a steep learning curve. We make a few final comments on barriers to adopting Lean from the perspective of a mathematics researcher, and suggestions on how researchers may begin to overcome them.

\begin{enumerate}
\item
For those who are not experienced in software development, installing Lean and its toolchain, keeping it updated to the latest version, and so on, can already be a nontrivial barrier. For researchers looking to try out Lean without installing it, we recommend the Lean playground \cite{playground}, a browser-based implementation of Lean with no installation requirements.

\item
\emph{Dependent type theory} (the basis for Lean's proof-checking logic) is still new to most researchers, and those who are not comfortable with programming will encounter significant barriers on this front too. Even for those who are comfortable with both, the syntax and naming conventions of Lean can be unclear for new users: the difference between \lstinline{true}, \lstinline{True} and \lstinline{isTrue} is opaque. A comprehensive, deep understanding of both is not strictly required to do mathematics in Lean, but they cannot be entirely glossed over.

There are many excellent community-written books aimed at Lean beginners with no knowledge of type theory, of which we mention only \cite{mil} and \cite{tpil}. The Lean Game Server \cite{gameserver}, a repository of interactive games aimed at teaching Lean syntax, can also serve as a good first introduction.

\item
Finding the ``best" definition of a concept requires familiarity not only with type theory and the syntax of Lean, but also with the library design of Lean and mathlib, as the user will certainly wish to lean on the available bank of already-proven lemmas. The available search tools are powerful, and the documentation is well-written, but it is worth spending time on these upfront rather than using them solely as a reference. mathlib also uses very systematic naming conventions for its lemmas, and learning these conventions can help significantly.

\item Lean 4 and mathlib4 are still in active development, and at time of writing, new versions are released once a month. Most changes add new functionality, and the occasional deprecation of functionality is minor, especially for new users: for instance, \lstinline{IsChain} (used above) is new as of September 2025, and the older proposition \lstinline{Chain} is now deprecated. This is unlikely to break old code (within reasonable timeframes), but keeping up with the latest versions of Lean 4 and mathlib4 will ensure that the user has access to the largest, most up-to-date bank of definitions and lemmas.

Black-box tactics such as \lstinline{simp} can be brittle, however, and proofs making heavy use of these may break as new versions of Lean are released. \lstinline{simp} is a useful tool, but rather a blunt instrument, and it is recommended that users learn to use \lstinline{simp only [}\textit{lemma\_{}names}\lstinline{]} for greater control.
\end{enumerate}

Another possible application of Lean we can see is in teaching, especially among students who have some basic familiarity with programming: an introductory course in logic or proofs may benefit from the high standards of formalisation and the immediate feedback from Lean's proof checker. Using mathlib, this could also draw on the students' previous knowledge of abstract algebra (proving that, for instance, $(ab)^{-1} = b^{-1}a^{-1}$ in a group is an easy but nontrivial task in Lean), elementary number theory or discrete mathematics. The book \cite{mil} is a good source of inspiration here.

\bibliographystyle{alpha}
\bibliography{main}

\end{document}